\documentclass{article}
\usepackage{epsfig}
\usepackage{amsmath}
\usepackage{amsfonts}

\newtheorem{theorem}{Theorem}[section]

\newtheorem{pr}[theorem]{Proposition}
\newtheorem{definition}[theorem]{Definition}
\newtheorem{example}[theorem]{Example}

\newtheorem{remark}[theorem]{Remark}

\begin{document}

\title{Statistical convergence of nets through directed sets}

\author{AR. Murugan$^1$, J. Dianavinnarasi $^2$ and C. Ganesa Moorthy $^3$\footnote {Corresponding author, e-mail: ganesamoorthyc@gmail.com} \ \\
%EndAName
$^{1,}$$^{2,}$$^{3}$ {\small Department of Mathematics, Alagappa University, Karaikudi-630 004, India.}\\}

\date{}
\maketitle

\begin{abstract}
The concept of statistical convergence based on asymptotic density is introduced in this article through nets. Some possible extensions of classical results for statistical convergence of sequences are obtained in this article, with extensions to nets.\\

$\textbf{Keywords.}$  Asymptotic density, Nets, Uniform spaces, Topological vector spaces.
\end{abstract}
\section{Introduction}
The concept of statistical convergence was introduced by H. Fast[7] and independently by H. Steinhaus in [23] as an applicable concept that generalizes the classical concept of usual convergence. This convergence was studied for sequences of numbers in [8, 9, 20], for sequences of elements in uniform spaces in[4, 16], for sequences of elements in paranormed spaces in [2, 11], for sequences of elements in topological groups in [6], for sequences of elements in metric spaces in [3], for sequences of elements in topological vector spaces in [15], and for sequences of elements in topological vector lattices in [1]. There are articles [17, 18] which study statistical convergence of double sequences and generalized sequences. There are generalizations of this concept through ideals in the articles [12, 13, 14, 21]. Almost all applicable statistical convergence ideas depend on asymptotic densities of sets. These sets may be subsets of $N,\;N\times N,\; N\times N\times N,\;...\;,$ where $N$ represents the set of all natural numbers. So, if the concept of asymptotic density for subsets of directed sets is introduced, then the concept of statistical convergence for nets can be introduced. This is done in the present article. For this purpose, a natural restriction is made on directed sets. The restriction is the following:
\par $\textit{For the directed sets $(D,\leq)$ considered in this article, to each $\beta\in D$,}$ $\textit{the set $\left\{\alpha\in D:\alpha\leq\beta\right\}$}$ \\$\textit{is finite and the set $\left\{\alpha\in D:\alpha\geq\beta\right\}$ is infinite}$.
\\It is assumed that all directed sets considered in this article satisfy this condition.
\par All directed sets considered through $N,\;N\times N,\;N\times N\times N,...,$ in earlier studies for statistical convergence satisfy this condition. Thus a common extension is proposed in this article.
\par There is an article [14] which discusses statistical convergence of nets through ideals, but not through a concept of asymptotic density. The present article presents statistical convergence of nets through a concept of asymptotic density for directed sets.
\par There are articles related to summability through statistical convergence (see [8, 22]) and articles for generalizations of asymptotic density (see [5]). Let us first introduce a concept of asymptotic density for our purpose.
\section{Asymptotic density}
\begin{definition}
Let $(D, \leq)$ be a directed set that satisfies the condition mentioned above.
\par To each $\alpha\in D$, let $D_{\alpha}=\{\beta\in D:\beta\leq\alpha\}$ and $|D_{\alpha}|$ denote the cardinality of $D_{\alpha}$. The lower asymptotic density of a nonempty subset $A$ of $D$ is defined as the number $\liminf\limits_{\alpha \in D}\frac{|A\cap D_{\alpha}|}{|D_{\alpha}|}$ and the upper asymptotic density of $A$ is defined as the number $\limsup\limits_{\alpha\in D}\frac{|A\cap D_{\alpha}|}{|D_{\alpha}|}$.
\par If the upper and lower densities are equal, then the common number is called the asymptotic density of $A$ and it is denoted by $\delta(A; D)$. Thus $\delta(A;D)=\lim\limits_{\alpha\in D}\frac{|A\cap D_{\alpha}|}{|D_{\alpha}|}$, in the real interval $[0,1]$. If $A$ is an empty subset, it is assumed that $\delta(A;D)=0$.
\par Here, for $x_{\alpha}\in \mathbb{R}$, the real line,
\begin{eqnarray*}
\liminf\limits_{\alpha\in D}x_{\alpha}&=&\sup\limits_{\beta\in D}\inf\limits_{\alpha\geq\beta}x_{\alpha}\;and\\
\limsup\limits_{\alpha\in D}x_{\alpha}&=&\inf\limits_{\beta\in D}\sup\limits_{\alpha\geq\beta}x_{\alpha}
\end{eqnarray*}
\end{definition}
\begin{example}
Let $D=\{(x_{1},x_{2},x_{3}):x_{i}\in N,\;i=1,2,3\}$. Define $\leq$ on $D$ by:\\$(x_{1}, x_{2}, x_{3})\leq (y_{1}, y_{2}, y_{3})$ if and only if $x_{1}\leq y_{1}$, $x_{2}\leq y_{2}$ and $x_{3}\leq y_{3}$. Then to each $(y_{1},y_{2},y_{3})\in D$, the set $\{(x_{1},x_{2},x_{3})\in D:(x_{1},x_{2},x_{3})\leq (y_{1},y_{2},y_{3})\}$ is finite, and it contains $y_{1}\cdot y_{2}\cdot y_{3}$ elements.\\Let $A=\{(x, x, x):x\in N\}$. Then $\delta(A;D)=0$.
\end{example}
\begin{example}
Let $D=N$ be the directed set with the usual order relation. Then to each $\alpha\in D$, $D_{\alpha}=\{\beta\in D: \beta\leq\alpha\}$ has precisely $\alpha$ elements. The asymptotic density introduced  in Definition 2.1 for $D$ coincides with the classical asymptotic density for subsets of $N$.
\end{example}
\begin{definition}
Let $(D^{(1)},\leq^{(1)})$ and $(D^{(2)},\leq^{(2)})$ be two directed sets. Let $D=D^{(1)}\times D^{(2)}$. Define the product order $\leq$ in $D$ by: $(x_{1},x_{2})\leq (y_{1},y_{2})$ if and only if $x_{1}\leq^{(1)}y_{1}$ and $x_{2}\leq^{(2)}y_{2}$. Observe again that to each $\alpha\in D$, the set $D_{\alpha}=\{\beta\in D:\beta\leq\alpha\}$ is finite. This definition can be extended to any  Cartesian product of a finite number of directed sets.
\end{definition}
\begin{remark}
If $A\subseteq D^{(1)}$, and if $\delta(A;D^{(1)})$ exists then $\delta(A\times D^{(2)};D)=\delta(A;D^{(1)})$; for the notations used in the previous Definition 2.4. Moreover, if $B\subseteq D^{(2)}$ and $\delta(A;D^{(1)})=0$ then $\delta(A\times B;D)=\delta(A;D^{(1)})$.
\end{remark}
\begin{pr}
Let $D$ be one among the directed sets $N, N\times N, N\times N\times N,\cdots ,$ when $N$ is endowed with the usual order, and the other sets are endowed with the corresponding product orders. Then, to each $\gamma\in D$, \\$\delta(\{\alpha\in D: \alpha\; not\;greater\;than\;or\;equal\;to\;\gamma\}; D)=0$, and hence $\delta(\{\alpha\in D:\alpha\geq\gamma\}; D)=1$.
\end{pr}
\textbf{Proof.}
It is easy to verify the relation $\delta(\{x\in D:\alpha\geq\gamma\};D)=1$.
\begin{example}
  Consider the set $N$ with the following different order relation. $m\leq n$ if and only if $m$ divides $n$. Then $N$ is a directed set with the properties mentioned in the introduction. Fix $k\in N\backslash\{1\}$. Let $A=\{n\in N: n \; is\; not\; greater\; than\; or\;equal\; to\; k\}=N\backslash\{k,\;2k,\;3k,\;\cdots\}=N\backslash kN$ (say).\\
  If $m\in N\backslash kN$ and if $D_{m}=\{n\in N:n\leq m\}$, then, for $i\in D_{m}$, $i\in N\backslash kN$ and hence $A\cap D_{m}=D_{m}$. This shows that $\limsup_{m\rightarrow\infty}\frac{|A\cap D_{m}|}{|D_{m}|}=1$.\\
  If $m=k^{i}$ for some $i>1$, then $D_{m}=\{n\in N:n\leq m\}=\{1,k,2k,\cdots,k^{i}\}$, when $A\cap D_{m}=\{1\}$. This shows that $\liminf_{m\rightarrow\infty}\frac{|A\cap D_{m}|}{|D_{m}|}=0$. In particular, $\delta(A;N)$ does not exist. However, if $A=\{n\in N:n\; is\; not \;greater\; than\;or\; equal\; to\; 1\}=N\backslash \{1\}$, then $\delta(A;N)=1$. Now, let $D=N\backslash \{1\}$ and consider the order relation defined above. To each $\beta\in D$, let $D_{\beta}=\{\alpha\in D:\alpha\leq \beta\}$. For a fixed $\gamma\in D$, let $B=\{\alpha\in D:\alpha\geq\gamma\}$. Then $\limsup_{\beta\in D}\frac{|B\cap D_{\beta}|}{|D_{\beta}|}=1>0$.
  \end{example}
  \begin{definition}
    A directed set is said to satisfy the condition (*), if to each fixed $\gamma\in D$, for the set $B=\{\alpha\in D:\alpha\geq \gamma\}$, it is true that $\limsup_{\beta\in D}\frac{|B\cap D_{\beta}|}{|D_{\beta}|}>0$, when $D_{\beta}=\{\alpha\in D:\alpha\leq\beta\}$.
  \end{definition}
\section{Statistical convergence}
\begin{definition}
Let $(x_{\alpha})_{\alpha\in D}$ be a net in a topological space $(X,\tau)$ and let $x\in X$. Let us say that $(x_{\alpha})_{\alpha\in D}$ converges statistically to $x$ in $(X,\tau)$, if, to each $U\in\tau$ such that $x\in U$, the relation $\delta(\{\alpha\in D:x_{\alpha}\notin U\};D)=0$ is true.
\end{definition}
\par Let us first verify the uniqueness of statistical limits in Hausdorff spaces.
\begin{pr}
Suppose $(x_{\alpha})_{\alpha\in D}$ be a net in a Hausdorff space $(X,\tau)$ such that it converges statistically to $x$ and $y$ in $X$. Then $x=y$.
\end{pr}
\textbf{Proof.}
Suppose $x\neq y$. Then there are disjoint open sets $U$ and $V$ such that $x\in U$ and $y\in V$. Then
\begin{eqnarray*}
 \{\alpha\in D:x_{\alpha}\notin U\}\cup \{\alpha\in D:x_{\alpha}\notin V\} &=& \{\alpha\in D:x_{\alpha}\notin U\cap V\}=D.
\end{eqnarray*}
But $\delta(\{\alpha\in D:x_{\alpha}\notin U\}\cup \{\alpha\in D:x_{\alpha}\notin V\};D)=0$ and $\delta(D;D)=1$; which is a contradiction. Therefore $x=y$. Observe that $\delta(A\cup B; D)=0$ whenever $\delta(A;D)=0$ and $\delta(B;D)=0$, for subsets $A$ and $B$ of $D$.

\begin{pr}
Let $(D^{(1)}, \leq^{(1)})$, $(D^{(2)}, \leq^{(2)})$ and $(D, \leq)$ be as in Definition 2.4. Let $(X, \tau_{X})$ and $(Y, \tau_{Y})$ be given topological spaces. Let $\tau$ be the product topology on $X\times Y$. Let $(x_{\alpha})_{\alpha\in D^{(1)}}$ and $(y_{\beta})_{\beta\in D^{(2)}}$ be two nets in $X$ and $Y$ respectively. Then $((x_{\alpha},y_{\beta}))_{(\alpha,\beta)\in D}$ converges statistically to some $(x,y)$ in $(X\times Y, \tau)$ if and only if $(x_{\alpha})_{\alpha\in D^{(1)}}$ converges statistically to $x$ in $(X, \tau_{X})$ and $(y_{\beta})_{\beta\in D^{(2)}}$ converges statistically to $y$ in $(Y, \tau_{Y})$.
\end{pr}
\textbf{Proof.}
  Suppose $(x_{\alpha})_{\alpha\in D^{(1)}}$ converges statistically to $x$ in $(X, \tau_{X})$ and $(y_{\beta})_{\beta\in D^{(2)}}$ converges statistically to $y$ in $(Y, \tau_{Y})$. Fix $U\in \tau_{X}$ and $V\in \tau_{Y}$ such that $x\in U$ and $y\in V$. Then
  \begin{eqnarray*}
    \delta(\{\alpha\in D^{(1)}:x_{\alpha}\notin U\};D^{(1)}) &=& 0\;\;and\\
    \delta(\{\beta\in D^{(2)}:y_{\beta}\notin V\};D^{(2)})  &=& 0.
  \end{eqnarray*}
By Remark 2.5,
\begin{eqnarray*}
  \delta(\{\alpha\in D^{(1)}:x_{\alpha}\notin U\}\times D^{(2)}\cup D^{(1)}\times\{\beta\in D^{(2)}:y_{\beta}\notin V\};D) &=& 0.
\end{eqnarray*}
Thus,
\begin{eqnarray*}
  \delta(\{(\alpha,\beta)\in D:(x_{\alpha},y_{\beta})\notin U\times V\};D) &=& 0.
\end{eqnarray*}
This implies that $((x_{\alpha},y_{\beta}))_{(\alpha,\beta)\in D}$ converges statistically to $(x,y)$ in $(X\times Y,\tau)$. Conversely, assume that $((x_{\alpha}, y_{\beta}))_{(\alpha, \beta)\in D}$ converges statistically to $(x, y)$ in $(X\times Y, \tau)$. Fix an open neighborhood $U$ of $x$ in $(X,\tau_{X})$. Then
\begin{eqnarray*}
  \delta(\{(\alpha,\beta)\in D:(x_{\alpha},y_{\beta})\notin U\times Y\};D) &=& 0.
\end{eqnarray*}
So, $\delta(\{\alpha\in D^{(1)}:x_{\alpha}\notin U\};D^{(1)}) = 0$. This implies that $(x_{\alpha})_{\alpha\in D^{(1)}}$ converges statistically to $x$ in $(X, \tau_{X})$. Similarly, $(y_{\beta})_{\beta\in D^{(2)}}$ converges statistically to $y$ in $(Y,\tau_{Y})$.
\begin{pr}
Let $(X,\tau_{X})$ , $(Y,\tau_{Y})$ and $(X\times Y, \tau)$ be as in the previous Proposition 3.3. Let $(x_{\alpha})_{\alpha\in D}$ be a net that converges statistically to some $x$ in $(X, \tau_{X})$, for some directed set $(D,\leq)$. Let $(y_{\alpha})_{\alpha\in D}$ be a net that converges statistically to some $y$ in $(Y,\tau_{Y})$. Then $((x_{\alpha}, y_{\alpha}))_{\alpha\in D}$ converges statistically to $(x,y)$ in $(X\times Y,\tau)$. On the other hand, if $((x_{\alpha},y_{\alpha}))_{\alpha\in D}$ converges statistically to some $(x,y)$ in $(X\times Y,\tau)$ then $(x_{\alpha})_{\alpha\in D}$ converges statistically to $x$ in $(X,\tau_{X})$ and $(y_{\alpha})_{\alpha\in D}$ converges statistically to $y$ in $(Y,\tau_{Y})$.
\end{pr}
\textbf{Proof.}
Suppose $(x_{\alpha})_{\alpha\in D}$ converges statistically to $x$ and $(y_{\alpha})_{\alpha\in D}$ converges statistically to $y$.\par Let $U$ be an open neighbourhood of $x$ in X and  $V$ be an open neighbourhood of $y$ in Y. Then
\begin{eqnarray*}
  \delta(\{\alpha\in D:x_{\alpha}\notin U\}\cup\{\alpha\in D:y_{\alpha}\notin V\};D) &=&0.
\end{eqnarray*}
\begin{eqnarray*}
That\;is,\; \delta(\{\alpha\in D:(x_{\alpha},y_{\alpha})\notin U\times V\};D) &=& 0.
\end{eqnarray*}
So, $((x_{\alpha},y_{\alpha}))_{\alpha\in D}$ converges statistically to $(x,y)$.\par Conversely assume that $((x_{\alpha},y_{\alpha}))_{\alpha\in D}$ converges statistically to $(x,y)$. Let $U$ be an open neighbourhood of $x$. Then
\begin{eqnarray*}
  \delta(\{\alpha\in D:(x_{\alpha},y_{\alpha})\notin U\times Y\};D)&=& 0.
\end{eqnarray*}
 That is, $\delta(\{\alpha\in D:x_{\alpha}\notin U\};D)=0$. This implies that $(x_{\alpha})_{\alpha\in D}$ converges statistically to $x$. Similarly $(y_{\alpha})_{\alpha\in D}$ converges statistically to $y$.
\begin{remark}
Proposition 3.3 and Proposition 3.4 can be extended to any  Cartesian product of a finite number of spaces.
\end{remark}
\begin{pr}
Let $(X,\tau_{X})$ and $(Y,\tau_{Y})$ be topological spaces and let $f:(X,\tau_{X})\rightarrow(Y,\tau_{Y})$ be a function which is continuous at a point $x$ in $X$. Let $(x_{\alpha})_{\alpha\in D}$ be a net that converges statistically to some $x$ in $(X,\tau_{X})$. Then $(f(x_{\alpha}))_{\alpha\in D}$ converges statistically to $f(x)$ in $(Y, \tau_{Y})$.
\end{pr}
\textbf{Proof.}
  Let $U$ be an open neighbourhood of $f(x)$ in $(Y,\tau_{Y})$. Then there is an open neighbourhood $V$ of $x$ in $(X,\tau_{X})$ such that $f(V)\subseteq U$. Then $\{\alpha\in D:f(x_{\alpha})\notin U\}\subseteq\{\alpha\in D:x_{\alpha}\notin V\}$ and $\delta(\{\alpha\in D:x_{\alpha}\notin V\};D)=0$. So, $\delta(\{\alpha\in D:f(x_{\alpha})\notin U\};D)=0$. This proves that $(f(x_{\alpha}))_{\alpha\in D}$ converges statistically to $f(x)$ in $(Y,\tau_{Y})$.
\begin{pr}
 Let $D^{(1)}, D^{(2)}$ and $D$ be as in Proposition 3.3. Let $(x_{\alpha})_{\alpha\in D^{(1)}}$ and $(y_{\beta})_{\beta\in D^{(2)}}$ be two nets in a topological vector space $X$ over the field of real numbers or the field of complex numbers. Let $(a_{\alpha})_{\alpha\in D^{(1)}}$ be a net of scalars. If $(x_{\alpha})_{\alpha\in D^{(1)}}$, $(y_{\beta})_{\beta\in D^{(2)}}$ and $(a_{\alpha})_{\alpha\in D^{(1)}}$ converge statistically to $x,y$ and $a$ respectively, then $(x_{\alpha}+y_{\beta})_{(\alpha,\beta)\in D}$ and $(a_{\alpha}y_{\beta})_{(\alpha,\beta)\in D}$ converge statistically to $x+y$ and $ay$ respectively.
 \end{pr}
\textbf{Proof.}
  Use Proposition 3.3 and Proposition 3.6. Observe that, it has been assumed that, the addition and the scalar multiplication in a topological vector space are jointly continuous.
\begin{pr}
Let $(x_{\alpha})_{\alpha\in D}$ and $(y_{\alpha})_{\alpha\in D}$ be two nets in a topological vector space $X$; with respect to a common directed set $D$. Let $(a_{\alpha})_{\alpha\in D}$ be a net of scalars. If $(x_{\alpha})_{\alpha\in D}$, $(y_{\alpha})_{\alpha\in D}$ and $(a_{\alpha})_{\alpha\in D}$ converge statistically to $x, y$ and $a$ respectively, then $(x_{\alpha}+y_{\alpha})_{\alpha\in D}$ and $(a_{\alpha}y_{\alpha})_{\alpha\in D}$ converge statistically to $x+y$ and $ay$ respectively.
\end{pr}
\textbf{Proof.}
  Use Proposition 3.4 and Proposition 3.6.
\begin{remark}
One may derive results similar to Proposition 3.7 and Proposition 3.8 for the structures, topological groups, topological rings, and topological algebras.
\end{remark}
\section{Statistically Cauchy nets}
The concept of statistically Cauchy nets is to be introduced for uniform spaces. For the concepts and notations in uniform spaces, one may refer to the book of Kelley [10] on General topology. The following definition agrees with the known definitions for statistically Cauchy sequences and statistically Cauchy double sequences(see [8, 17, 19]).
\begin{definition}
Let $(X,\mathfrak{U})$ be a uniform space with a uniformity $\mathfrak{U}$. A net $(x_{\alpha})_{\alpha\in D}$ in X is said to be statistically Cauchy if, for given $U\in \mathfrak{U}$, there is a $\gamma\in D$ such that
\begin{eqnarray*}
  \delta(\{\alpha\in D:(x_{\alpha}, x_{\gamma})\notin U, \alpha\geq\gamma\};D) &=& 0.
\end{eqnarray*}
\end{definition}
\par It is easy to verify that every Cauchy net is a statistically Cauchy net, and hence every converging net is a statistically Cauchy net in a uniform space. It is also possible to prove that statistical convergence implies statistical Cauchyness in a uniform space.
\begin{pr}
Let $D$ be a directed set. Then every statistically convergent net $(x_{\alpha})_{\alpha\in D}$ in a uniform space is statistically Cauchy.
\end{pr}
\textbf{Proof.}
  Let $(x_{\alpha})_{\alpha\in D}$ be a net which converges statistically to $x$ in a uniform space $(X,\mathfrak{U})$. Fix $U\in \mathfrak{U}$. Find a symmetric $V\in \mathfrak{U}$ such that $V\circ V\subseteq U$. For this $V$, $\delta(\{\alpha\in D:(x_{\alpha}, x)\notin V\}; D)=0$ and hence there is a $\gamma\in D$ such that $(x_{\gamma}, x)\in V$. Then $\{\alpha\in D: (x_{\alpha}, x_{\gamma})\notin U\}\subseteq \{\alpha\in D:(x_{\alpha}, x)\notin V\}$. Thus, $\delta(\{\alpha\in D:(x_{\alpha},x_{\gamma})\notin U, \alpha\geq\gamma\};D)=0$. This proves that $(x_{\alpha})_{\alpha\in D}$ is statistically Cauchy.\\

 Let us recall the order in product of two directed sets  described in Definition 2.4.
\begin{pr}
Let $(x_{\alpha})_{\alpha\in D}$ be a net that is statistically Cauchy in a uniform space $(X,\mathfrak{U})$. Then for given $U\in \mathfrak{U}$, there is a $\gamma\in D$ such that\\$\delta(\{(\alpha,\beta)\in D\times D:(x_{\alpha},x_{\beta})\notin U, \alpha\geq\gamma, \beta\geq\gamma\};D\times D)=0$.
\end{pr}
\textbf{Proof.}
  Fix $U\in\mathfrak{U}$. Find a symmetric $V\in\mathfrak{U}$ such that $V\circ V\subseteq U$. For this V, there is a $\gamma\in D$ such that $\delta(\{\alpha\in D:(x_{\alpha},x_{\gamma})\notin V,\alpha\geq\gamma\};D)=0$. Since
  \begin{eqnarray*}
    \{(\alpha,\beta)\in D\times D:(x_{\alpha},x_{\beta})\notin U,\alpha\geq\gamma,\beta\geq\gamma\} &\subseteq & \{(\alpha,\beta)\in D\times D:(x_{\alpha},x_{\gamma})\notin V \\
     &&\quad or\;(x_{\beta},x_{\gamma})\notin V, \alpha\geq\gamma, \beta\geq\gamma\} \\
     &\subseteq& (\{\alpha\in D:(x_{\alpha}, x_{\gamma})\notin V, \alpha\geq\gamma\}\times D)\\
      &&\quad \cup (D\times \{\beta\in D:(x_{\beta},x_{\gamma})\notin V, \beta\geq\gamma\}),
     \end{eqnarray*}
      by Remark 2.5,
      \begin{eqnarray*}
        \delta(\{(\alpha,\beta)\in D\times D:(x_{\alpha},x_{\beta})\notin U, \alpha\geq\gamma, \beta\geq\gamma\};D\times D) &=& 0.
      \end{eqnarray*}
\begin{pr}
Let $D^{(1)},D^{(2)}$ and $D$ be as in Proposition 3.3. Let $(X,\mathfrak{U}_{X})$ and $(Y,\mathfrak{U}_{Y})$ be two uniform spaces. Let $\mathfrak{U}$ be the product uniformity on $X\times Y$. Let $(x_{\alpha})_{\alpha\in D^{(1)}}$ and $(y_{\beta})_{\beta\in D^{(2)}}$ be two nets in $X$ and $Y$ respectively. Then $((x_{\alpha},y_{\beta}))_{(\alpha,\beta)\in D}$ is statistically Cauchy in $(X\times Y,\mathfrak{U})$ if $(x_{\alpha})_{\alpha\in D^{(1)}}$ is statistically Cauchy in $(X, \mathfrak{U}_{X})$ and $(y_{\beta})_{\beta\in D^{(2)}}$ is statistically Cauchy in  $(Y, \mathfrak{U}_{Y})$. Moreover, if $D^{(1)}$ and $D^{(2)}$ satisfy the condition (*) mentioned in Definition 2.8, and $((x_{\alpha}, y_{\beta}))_{(\alpha, \beta)\in D}$ is statistically Cauchy in $(X\times Y, \mathfrak{U})$, then $(x_{\alpha})_{\alpha\in D^{(1)}}$ is statistically Cauchy in $(X, \mathfrak{U}_{X})$ and $(y_{\beta})_{\beta\in D^{(2)}}$ is statistically Cauchy in  $(Y, \mathfrak{U}_{Y})$.
\end{pr}
\textbf{Proof.}
  The proof follows from the set relation: For $U\in \mathfrak{U}_{X}$, $V\in \mathfrak{U}_{Y}$, $\gamma_{1}\in D^{(1)}$ and for $\gamma_{2}\in D^{(2)}$, it is true that
  \begin{eqnarray*}
  \{(\alpha, \beta)\in D: ((x_{\alpha}, x_{\gamma_{1}}), (y_{\beta}, y_{\gamma_{2}}))&\notin& U\times V, (\alpha, \beta)\geq(\gamma_{1},\gamma_{2})\}\\
  &=& (\{\alpha\in D^{(1)}:(x_{\alpha}, x_{\gamma_{1}})\notin U, \alpha\geq\gamma_{1}\}\times \{\beta\in D^{(2)}:\beta\geq\gamma_{2}\})\\
  &&\quad \cup (\{\alpha\in D^{(1)}: \alpha\geq \gamma_{1}\}\times\{\beta\in D^{(2)}:(y_{\beta}, y_{\gamma_{2}})\notin V, \beta\geq \gamma_{2}\}).
  \end{eqnarray*}

\begin{pr}
Let $(X,\mathfrak{U}_{X})$ and $(Y,\mathfrak{U}_{Y})$ be two uniform spaces. Let $\mathfrak{U}$ be the product uniformity on $X\times Y$. Let $(x_{\alpha})_{\alpha\in D}$ and $(y_{\alpha})_{\alpha\in D}$ be nets in X and Y respectively. Then $((x_{\alpha},y_{\alpha}))_{\alpha\in D}$ is statistically Cauchy in $X\times Y$ if and only if $(x_{\alpha})_{\alpha\in D}$ is statistically Cauchy in X and $(y_{\alpha})_{\alpha\in D}$ is statistically Cauchy in Y.
\end{pr}
\textbf{Proof.}
Suppose $(x_{\alpha})_{\alpha\in D}$ and $(y_{\alpha})_{\alpha\in D}$ be statistically Cauchy. Fix $U\in \mathfrak{U}_{X}$ and $V\in\mathfrak{U}_{Y}$. Then there is a $\gamma\in D$ such that $\delta(\{\alpha\in D:(x_{\alpha},x_{\gamma})\notin U,\alpha\geq\gamma\};D)=0$ and $\delta(\{\alpha\in D:(y_{\alpha},y_{\gamma})\notin V,\alpha\geq\gamma\};D)=0$. The statistically Cauchyness  of $((x_{\alpha},y_{\alpha}))_{\alpha\in D}$ follows from the relation:\\
\begin{eqnarray*}
\{\alpha\in D:(x_{\alpha},x_{\gamma})\notin U&or&(y_{\alpha}, y_{\gamma})\notin V,\;\alpha\geq \gamma\} \\
&&\quad\subseteq \{\alpha\in D: (x_{\alpha},x_{\gamma})\notin U, \alpha\geq \gamma\}\cup\{\alpha\in D: (y_{\alpha},y_{\gamma})\notin V, \alpha\geq\gamma\}.
\end{eqnarray*}
\par Conversely, assume that $((x_{\alpha}, y_{\alpha}))_{\alpha\in D}$ is statistically Cauchy. Fix $U\in\mathfrak{U}$. Then there is a $\gamma\in D$ such that
\begin{eqnarray*}
  \delta(\{\alpha\in D: (x_{\alpha},x_{\gamma})\notin U, \alpha\geq\gamma\};D) =\delta(\{\alpha\in D: (x_{\alpha},x_{\gamma})\notin U\;or\;(y_{\alpha},y_{\gamma})\notin Y\times Y,\alpha\geq\gamma\}; D)=0.
\end{eqnarray*}
This shows that $(x_{\alpha})_{\alpha\in D}$ is statistically Cauchy. Similarly $(y_{\alpha})_{\alpha\in D}$ is statistically Cauchy.
\begin{pr}
Let $f:(X, \mathfrak{U})\rightarrow(Y, \mathfrak{V})$ be a uniformly continuous function from a uniform space $(X, \mathfrak{U})$ into a uniform space $(Y, \mathfrak{V})$. Let $(x_{\alpha})_{\alpha\in D}$ be a statistically Cauchy net in $(X, \mathfrak{U})$. Then $(f(x_{\alpha}))_{\alpha \in D}$ is a statistically Cauchy net in $(Y, \mathfrak{V})$.
\end{pr}
\textbf{Proof.}
  Fix $V\in\mathfrak{V}$. Find a $U\in\mathfrak{U}$ such that $(f(x),f(y))\in V$, whenever $(x, y)\in U$. Find a $\gamma\in D$ such that $\delta(\{\alpha\in D: (x_{\alpha},x_{\gamma})\notin U, \alpha\geq\gamma\};D)=0$. Then $\delta(\{\alpha\in D: (f(x_{\alpha}),f(x_{\gamma}))\notin V, \alpha\geq\gamma\};D)=0$, because $\{\alpha\in D: (f(x_{\alpha}),f(x_{\gamma}))\notin V, \alpha\geq\gamma\}\subseteq\{\alpha\in D: (x_{\alpha},x_{\gamma})\notin U, \alpha\geq\gamma\}$.
\begin{remark}
Let $(X, \tau)$ be a topological vector space. The usual uniformity on $X$ implies the following: A net $(x_{\alpha})_{\alpha\in D}$ is Cauchy in $X$ if and only if for every neighbourhood $U$ of $0$ there is a $\gamma\in D$ such that
\begin{eqnarray*}
  \delta(\{\alpha\in D: x_{\alpha}-x_{\gamma}\notin U, \alpha\geq \gamma\}; D) &=& 0.
\end{eqnarray*}
One can derive the following Proposition 4.8 and Proposition 4.9 which are similar to Proposition 3.7 and Proposition 3.8.
\end{remark}
\begin{pr}
Let $D^{(1)}, D^{(2)}$, $D$, $(x_{\alpha})_{\alpha\in D^{(1)}},\;(y_{\beta})_{\beta\in D^{(2)}},\;(a_{\alpha})_{\alpha\in D^{(1)}} $ and $X$ be as in Proposition 3.7. Let $x\in X$ and $a$ be a scalar. If $(x_{\alpha})_{\alpha\in D^{(1)}},\;(y_{\beta})_{\beta\in D^{(2)}}\;and\;(a_{\alpha})_{\alpha\in D^{(1)}}$ are statistically Cauchy, then $(x_{\alpha}+y_{\beta})_{(\alpha, \beta)\in D}$, $(a_{\alpha}x)_{\alpha\in D^{(1)}}$ and $(ax_{\alpha})_{\alpha\in D^{(1)}}$ are statistically Cauchy.
\end{pr}
\textbf{Proof.}
  Use Proposition 4.6 and Proposition 4.4.
\begin{pr}
Let $(x_{\alpha})_{\alpha\in D}$, $(y_{\alpha})_{\alpha\in D}$ and $X$ be as in Proposition 3.8. If $(x_{\alpha})_{\alpha\in D}$ and $(y_{\alpha})_{\alpha\in D}$ are statistically Cauchy, then $(x_{\alpha}+y_{\alpha})_{\alpha\in D}$ is statistically Cauchy.
\end{pr}
\textbf{Proof.}
  Use Proposition 4.5 and Proposition 4.6.
\section{Net Spaces}
Corresponding to sequence spaces, net spaces can be constructed. The following construction is similar to the construction given in [20].
The following construction uses the Propositions 3.7, 3.8, 4.8 and 4.9. Since verifications part is a direct one, it is omitted.
\par Let $(X, \tau)$ be a topological vector space with the natural uniformity $\mathfrak{U}$ that induces the topology $\tau$. Let $D$ be a fixed directed set. Let $M = \{(x_{\alpha})_{\alpha\in D}: \{x_{\alpha}: \alpha\in D\}is\;a\;bounded\;subset\;of\;X\}$.
\begin{eqnarray*}
  Let\;M_{cy} &=& \{(x_{\alpha})_{\alpha\in D}\in M:\;(x_{\alpha})_{\alpha\in D}\;is\;statistically\;Cauchy\}. \\
  Let\;M_{ct} &=& \{(x_{\alpha})_{\alpha\in D}\in M:\;(x_{\alpha})_{\alpha\in D}\;converges\;statistically\;in\;X\}. \\
 Let\;M_{0} &=& \{(x_{\alpha})_{\alpha\in D}\in M:\;(x_{\alpha})_{\alpha\in D}\;converges\;statistically\;to\;zero\;in\;X\}.
\end{eqnarray*}
To each balanced neighbourhood $U$ of zero in $X$, define a function $p_{U}$ on $M$ by
\begin{eqnarray*}
  p_{U}((x_{\alpha})_{\alpha\in D}) &=& \sup \{\lambda\geq 0: \lambda x_{\alpha}\in U, \forall\alpha\in D\},
\end{eqnarray*}
and define a subset $N_{U}$ of $M$ by
\begin{eqnarray*}
  N_{U} &=& \{(x_{\alpha})_{\alpha\in D}\in M: p_{U}((x_{\alpha})_{\alpha\in D})< 1\}.
\end{eqnarray*}
Then the collection of the sets of the form $N_{U}$ forms a local base for $M$ that makes $M$ into a topological vector space under pointwise algebraic operations. Also $M_{cy}$ is a closed linear subspace of $M$. If $(X, \mathfrak{U})$ is a complete topological vector space, then $M$ is a complete topological vector space and $M_{ct}$ and $M_{0}$ are closed linear subspaces of $M$.

\end{document}